\documentclass{article}
\usepackage{amsmath}
\usepackage{mathrsfs}
\usepackage{amsfonts}
\usepackage{amssymb}
\usepackage{graphicx}
\usepackage{amsthm}
\usepackage[toc,page]{appendix}
\usepackage{geometry}
\usepackage[all]{xy}
\usepackage{lmodern}
\usepackage[T1]{fontenc}
\usepackage[dvipsnames]{xcolor}
\usepackage{graphics}
\usepackage{comment}
\usepackage{csquotes}
\usepackage[colorlinks=true,linkcolor=NavyBlue,urlcolor=RoyalBlue,citecolor=PineGreen,linktocpage=true]{hyperref}
\usepackage{enumerate}
\usepackage{caption}
\usepackage{subcaption}
\usepackage{float}
\usepackage{multicol}
\usepackage{algpseudocodex}
\usepackage{algorithm}
\usepackage{tcolorbox}

\numberwithin{equation}{section}

\newtheorem{Theorem}{Theorem}[section]	
\newtheorem*{Theorem*}{Theorem}
\newtheorem{proposition}[Theorem]{Proposition} 
\newtheorem{lemma}[Theorem]{Lemma}
	
\newtheorem{defn}[Theorem]{Definition}
\newtheorem*{defn*}{Definition}

\newtheorem*{remark}{Remark}

\newtheorem{example}{Example}[section]

\usepackage{tikz}
\usetikzlibrary{decorations,arrows}
\usetikzlibrary{decorations.pathmorphing}
\usepgflibrary{decorations.pathreplacing} 
\usepackage{pgf}
\usetikzlibrary{patterns}

\tikzset{
	schraffiert/.style={pattern=horizontal lines,pattern color=#1},
	schraffiert/.default=black
}

\tikzset{
	ultra thin/.style= {line width=0.1pt},
	very thin/.style=  {line width=0.2pt},
	thin/.style=       {line width=0.4pt},
	semithick/.style=  {line width=0.6pt},
	thick/.style=      {line width=0.8pt},
	very thick/.style= {line width=1.2pt},
	ultra thick/.style={line width=2.4pt}
}

\usetikzlibrary{shapes}
\usetikzlibrary{positioning}

\title{Approximations of symbolic substitution systems in one dimension}
\author{Lior Tenenbaum}

\begin{document}
	\date{}
	\maketitle
	 \begin{abstract}
		Periodic approximations of quasicrystals are a powerful tool in analyzing spectra of Schrödinger operators arising from quasicrystals, given the known theory for periodic crystals.
		Namely, we seek periodic operators whose spectra approximate the spectrum of the limiting operator (of the quasicrystal). This naturally leads to study the convergence of the underlying dynamical systems.\\
		We treat dynamical systems which are based on one-dimensional substitutions. We first find natural candidates of dynamical subsystems to approximate the substitution dynamical system. Subsequently, we offer a characterization of their convergence and provide estimates for the rate of convergence.
		We apply the proposed theory to some guiding examples.
	\end{abstract}
	\section{Introduction}\label{Intro-Sect}
	In studying a Schrödinger operator coming from a quasicrystal, one often turns to study it by periodic crystals approximating  the underlying dynamical structure. Examples of this approach can be seen in earlier works such as \cite{Ostlund_1985, MDuneau_1989, C.Sire_1990, Kazou91,TINKAGAMMEL93}, and more recently in \cite{Jagannathan08, Akker14, Jagannathan15, 2DInsulator, Colbrook19, BBdN20}. This is done using finite volume approximants with either open, periodic or twisted boundary conditions, while trying to minimze the effects of the boundary conditions. In this paper, we deal with infinite approximants having periodic potential used to estimate a  Schrödinger operator coming from an aperiodic configuration of atoms on the infinite lattice $\mathbb{Z}$. These infinite periodic approximants are relatively well understood using Bloch-Floquet theory, which allows us to study them via finite volume operators with twisted boundary conditions. See for example \cite{BlochReview} or \cite{TwistLattice05}. The Schrödinger operators  we consider, which are simple cases of the tight binding model, are given by
	\begin{equation} \label{Schro-defn}
		H_\omega = -\Delta + V_\omega,
	\end{equation}
	where $\Delta$ is the discrete Laplacian and $V_\omega$ is a diagonal operator whose diagonal elements come from the underlying configuration. These 
	Recent results by Beckus, Bellissard, Cornean and Takase, \cite[Propostion 1.1]{BBC19} and \cite[Theorem 2.1]{BecTak21}, ensure that the distance between spectra coming from two underlying configurations are at most proportional to the distance between the hulls of the configurations. For a precise version of these results see Proposition \ref{Lipschitz-Bound} in the appendix.
	Given the aforementioned spectral approximation results, it is therefore natural to ask what are the restrictions on approximating aperiodic systems from periodic ones? And at which rate?\\
	We discuss in this paper such problems for one-dimensional aperiodic structures with substitution symmetry, also called inflation symmetry in other sources, see  \cite{MDuneau_1989}.  We restrict our attention to materials coming from an aperiodic configuration of atoms on a lattice by a substitution. We encode the atoms in the configurations by a finite set of symbols, called an alphabet and denoted by $\mathcal{A}$. The distance between the hulls we consider is inverse to the scale of local patterns on which they agree, and we denote this distance by $d_\mathcal{H}$. A precise definition of $d_\mathcal{H}$ and its characterization are discussed in Section \ref{dynam-sect} in the appendix. The local patterns in our hulls are finite words over $\mathcal{A}$. Explicit details and general theory for this setting can be found in \cite{DamFil}. \\
	In this setting, and using the results in \cite{BBC19} and  \cite{BecTak21}, we consider an approximation scheme for substitution aperiodic structures, using only the substitution structure similar to \cite{BBdN20}. Namely, we estimate an aperiodic structure by taking a periodic configuration $\omega_0$ and iteratively applying the substitution to $\omega_0$. This yields an iterative hull sequence (\textbf{IHS}), which we use to approximate the hull of a  substitution aperiodic structure. The hull of a configuration with substitution symmetry, is called a substitution subshift and denoted $\Omega(S)$. In this paper we discuss the IHS approximation scheme and the rates at which the scheme approaches $\Omega(S)$. We restrict our attention to the class of primitive substitutions, due to their nice dynamical properties. Explicit definitions and outlines of proofs are given in the appendix.\\
	We emphasize that this paper is part of a larger collaborative work with Ram Band, Siegfried Beckus and Felix Pogorzelski. This collaboration deals with similar issues in a generalized setting, of Cayley graphs of lattices in homogeneous Lie groups, and will appear in \cite{Curr}. 
	The results presented here deal with the one-dimensional case,  which do not fall under the framework of block substitutions used in \cite{Symbolic} and \cite{Curr}.
	\section{The approximation scheme} \label{scheme-sect}
\begin{figure}[H]
			\begin{center}
				\vspace{1.1cm}
				\includegraphics[width=0.5\linewidth]{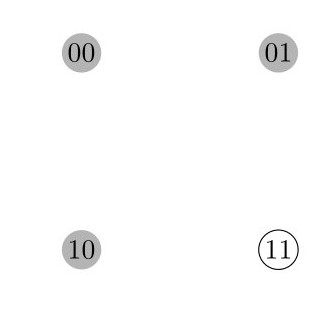}
				\vspace{0.75cm}
				\caption{\centering$G(S_{Fib})$ for the Fibonacci and the period doubling substitutions.}
				\label{Fibo-graph}
			\end{center}
		\end{figure}

\begin{figure}[H]
			\begin{center}
				\includegraphics[width=0.8\linewidth]{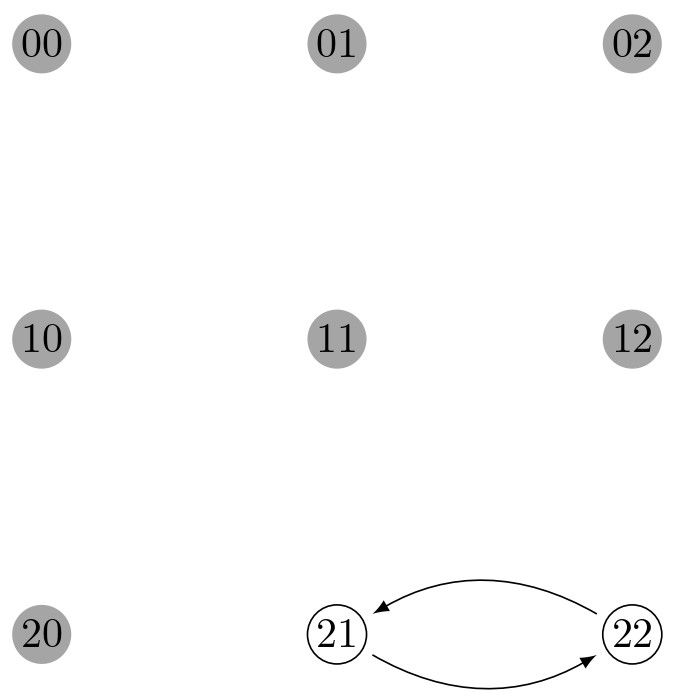}
				\caption{ \centering$G(S_{CE})$ for the substitution coming from the rule in \eqref{Poniros-Rule}.}
				\label{Poniros-graph} 
			\end{center}
			
		\end{figure}

	Two examples for substitutions, which we use to illustrate our theory, are the Fibonacci substitution given by the rule
	\begin{equation} \label{Fibo-Rule}
		0 \mapsto 01 \quad \text{and} \quad 1\mapsto 0,
	\end{equation}
	and 
	a counterexample substitution given by the rule
	\begin{equation} \label{Poniros-Rule}
		0\mapsto 001, \quad 1 \mapsto 200 \quad \text{and} \quad 2 \mapsto 102.   
	\end{equation}	
	Given the characterization of symbolic hulls distance via their local patterns, see for example \cite[Proposition 1.3.4]{LindMarc}, we seek to understand how these patterns change through our IHS approximations.  The local patterns for a substitution subshift are the ones that occur naturally through iterations of the substitution. For example, for the Fibonacci substitution given in \eqref{Fibo-Rule}, one can obtain its local patterns using the following derivation chain iterating the substitution.
	\begin{equation} \label{deriv-chain}
		1\; \mapsto\; 0 \mapsto\; 01 \; \mapsto \; 010 \; \mapsto\; 01001 \; \mapsto \; 01001010\; \mapsto \; 01 0 01 01 0 01 0 01 \; \mapsto \;  ...
	\end{equation}
	Local patterns occurring as a subword in this chain, are the local patterns appearing in $\Omega(S_{Fib})$. They are also called \textbf{legal} words or patterns. One can verify that the legal $2$-words in the Fiboncacci substitution case are $W(S_{Fib})_2=\{ 00,01,10 \}$ and the defective $2$-word is $11$. Relying on \cite[Corollary 6.2.5]{BeckusThesis} and \cite[Proposition 2.11]{BBdN20}, we know that the IHS  $\Omega_n(\omega_0)$ coming from a starting configuration $\omega_0$, converge if the $2$-words in $\omega_0$ are legal. Therefore, taking the constant configuration $\omega_0= 0^\infty$, yields a sequence of periodic hulls converging to $\Omega(S_{Fib})$.\\
	Building on this criterion, if the sequence of hulls eventually does not contain any defective $2$-word, then the IHS converge. For example, apply the substitution given in \eqref{Poniros-Rule}. Using a similar derivation chain to the one in \eqref{deriv-chain}, we conclude in this case that $W(S_{CE})_2=\{ 00,01,02,10,11,12,20 \}$, and the defective length $2$-words are $\{21,22\}$.  Starting with the periodic configuration $\omega_0=2^\infty$ and applying the substitution iteratively, one can see that the defective length $2$-words appears infinitely many times in the IHS. Therefore, the hulls coming from $\omega_0=2^\infty$ are bad approximants for $\Omega(S_{CE})$.
	To summarize, our approximation scheme is as follows.
	\begin{tcolorbox}
		\begin{center}
			{\Large The IHS approximation scheme} 
		\end{center}
	\begin{itemize}
		 \item \underline{Needed:}  A primitive substitution $S$ and an initial configuration $\omega_0$ over the lattice $\mathbb{Z}$.
		\item \underline{Iterative step:} At the $n$-th step, consider the symbolic hull $\Omega_n(\omega_0)$ of $S^n(\omega_0)$, obtained from applying $S$ $n$-times  to $\omega_0$.
		\item \underline{Result:} Determine how does the distance between the IHS and $\Omega(S)$ changes.
	\end{itemize}	 
	\end{tcolorbox}
	Following this discussion, we generate a directed graph, denoted $G(S)$, to keep track of defective $2$-words occurring in the hulls along the iteration sequence.
	The vertices of such a $G(S)$ graph are all the $2$-words over $\mathcal{A}$. We then draw edges between a defective word $u$ to a defective word $v$ occurring in the substitution of $S(u)$ to pay attention to the defective $2$-words in the hulls.
	Drawing these graphs for our examples, we obtain Figures \ref{Fibo-graph} and \ref{Poniros-graph}, with the legal $2$-words colored in gray for each substitution.
	\begin{multicols}{2}
		
		\columnbreak

	\end{multicols}
	Following \cite[Proposition 2.11]{BBdN20}, and using the proposed graphs, we obtain the following algorithmic result to determine whether our approximation scheme succeeds.
	\begin{Theorem} \label{1D-necess&suff}
		Let $S$ be a primitive one-dimensional substitution and $\omega_0$ be some initial configuration over $\mathbb{Z}$.  The following conditions are equivalent.
		\begin{enumerate}[(i)]
			\item\hspace{-0.2cm} \label{1D-conv} The IHS $\Omega_n(\omega_0)$ converge to $\Omega(S)$.
			\item \hspace{-0.2cm} \label{1D-finite-path} Any directed path in $G(S)$ starting in a $2$-word of $\omega_0$ does not contain a closed subpath.
			\item \hspace{-0.2cm} \label{1D-path-length} Any directed path in $G(S)$ starting in a $2$-word of $\omega_0$ has length strictly less than $\vert \mathcal{A} \vert^2$.
		\end{enumerate}
	\end{Theorem}
	The interested reader can find a proof of Theorem \ref{1D-necess&suff} in the appendix. 
	We say that an initial configuration $\omega_0$ and the IHS coming from it are \textbf{good} approximants if they satisfy the conditions in Theorem \ref{1D-necess&suff}, and \textbf{bad} approximants otherwise.
	Turning back to our motivating examples, we can see that any starting configuration is good in the case of the Fibonacci substitution. More generally, any starting configuration in the case of a self-correcting substitution, see \cite{Gahler} for an explicit definition, is a good approximant and yields a good approximating IHS. The class of self-correcting substitutions constitute the majority of one-dimensional substitution examples in the literature, and include the Fibonacci, Thue--Morse, Period--doubling and the Golay--Rudin--Shapiro substitutions. See \cite[Chapter 4]{BaakeGrimm13} for more details.
	On the other hand, considering the counterexample from \eqref{Poniros-Rule}, we see that an initial configuration, and its derived sequence of hulls, is bad if and only if it contains a defective $2$-word $21$ or $22$. More generally, any starting configuration containing a defective $2$-word in the case of a marked substitution, see \cite{Frid99} for a definition, will generate bad approximating IHS.
	\section{Approximation rates from the scheme} \label{Approx-rate-sect}
	We now turn to estimating the rates at which the approximation scheme gives, when the IHS approximants are good. Note that the IHS $\Omega_n(\omega_0)$ are bad if and only if $d_{\mathcal{H}}\big(  \Omega_n(\omega_0), \Omega(S) \big)$ is bounded away from $0$. To that end we need to consider how fast does the IHS agree with the substitution subshift on local patterns. Recall that a substitution $S$ has a matrix $M_S$, consult \cite[Section 5.3]{Queff10}, encoding the substitutions rule. The matrices for the substitutions given in \eqref{Fibo-Rule} and \eqref{Poniros-Rule} are
	\begin{equation}\label{Sub-Mat}
		M_{S_{Fib}}= \begin{pmatrix}
			1 && 1 \\
			1 && 0
		\end{pmatrix} \quad \text{and} \quad M_{S_{CE}}=\begin{pmatrix}
		2 && 2 && 1 \\
		1 && 0 && 1 \\
		0 && 1 && 1 
		\end{pmatrix}
	\end{equation}
	accordingly.
	A substitution being primitive, implies that the associated substitution matrix $M_S$ has a leading simple eigenvalue, using Perron--Frobenius (\textbf{PF}) theory, consult \cite[Chapter 8]{MeyerBook} or \cite[Section 5.3]{Queff10} for more details. This leading eigenvalue is called the PF-eigenvalue, which we denote $\theta_S$. The PF-eigenvalue for the matrices in \eqref{Sub-Mat} are $\frac{1+\sqrt{5}}{2}$ and $3$ accordingly.   The PF-theory implies that there exist constants $\hat{C}(S)>0$ and $\check{C}(S)>0$ such that 
	\begin{equation} \label{Perron-estim}
		\check{C}(S) \theta_S^n \leq  \vert S^n(a)\vert \leq \hat{C}(S)\theta_S^n \quad \text{for all} \quad a\in \mathcal{A} \quad \text{and}  \quad n\in \mathbb{N},
	\end{equation}
	where $\vert u\vert$ is the length of the word $u$.
	It is well known, see \cite{Sol98} or \cite{Dur00}, that a primitive substitution induces a linearly repetitive subshift. i.e., there exists a constant $C_S\geq 1$ for $\Omega(S)$, such that any legal word of length $\ell$ occurs inside all legal words of length greater than $\ell C_S$. We can show, see Section \ref{Rate-Sect} in the appendix, that the IHS $\Omega_n(\omega_0)$ agree on $\ell$-length  words when $n\geq \frac{\log \ell}{\log \theta_S}+ \tilde{C}_1$, for some appropriate constant $\tilde{C_1}>0$.
	This establishes an upper bound for the rate of convergence in \eqref{Conv-rate}.\\ 
	Conversely, when the initial configuration is periodic, a lower bound on the rate of the convergence follows from the number of $\ell$-length words occurring in $\Omega(S)$. Using a result from Coven--Hedlund \cite{Coven73}, we see that the number of $\ell$-words in $\Omega(S)$ must be at least $\ell+1$. One can see that a periodic hull can have, at any fixed word length $\ell$, at most as many different words as its period size. When $\omega_0$ is a periodic configuration and $\Omega_n(\omega_0)$ agrees with $\Omega(S)$ on $\ell$-length words, it follows that $n\geq \frac{\log \ell}{\log \theta_S}+ \tilde{C}_2(\omega_0;S)$, for an appropriate constant $\tilde{C}_2(\omega_0;S)>0$. Combining the bounds on $n$, we conclude that for a good periodic starting configuration, we get that the IHS generated from the approximation scheme satisfy
	\begin{equation} \label{Conv-rate}
		  \frac{C_2(\omega_0;S)}{\theta_S^n} \leq  d_{\mathcal{H}}\big( \Omega_n(\omega_0),\Omega(S) \big) \leq \frac{C_1(S)}{\theta_S^n} \quad \text{for all} \quad n\in \mathbb{N},
	\end{equation} 
	where $C_1(S)>0$ and $C_2(\omega_0;S)>0$ are appropriate constants. In particular, any periodic good approximating IHS have an asymptotically optimal distance from $\Omega(S)$ of $\theta_S^{-n}$.\\
	Applying the estimates in \eqref{Conv-rate} to the Fibonacci substitution gives us that any IHS converge to $\Omega(S_{Fib})$ in an asymptotic rate of $\big( \frac{1+\sqrt{5}}{2} \big)^{-n}$. More generally, for any self-correcting substitutions with PF-eigenvalue $\theta_S$ and any initial configuration $\omega_0$, we have that $d_\mathcal{H}\big( \Omega_n(\omega_0),\Omega(S) \big)$ is at most proportional to $\theta_S^{-n}$.\\
	Finally, turning back to the spectrum of a Schrödinger operator described in \eqref{Schro-defn}, the sequence of spectra arising from a good approximating IHS are at distance at most proportional to  $\theta_S^{-n}$. This shows that spectra of good approximating IHS approach the spectrum of a Schrödinger operator, with substitution symmetric potential, exponentially fast. In particular, for any self-correcting substitution and any initial configuration, the generated IHS yield an exponentially fast estimate for the spectrum of a Schrödinger operator.  A precise phrasing  of this statement is given in  Proposition \ref{Approx-Method} in the appendix.
	This can be used to prove the existence of spectral gaps similar to arguments in \cite{Hege22}.
	
	\subsection*{Acknowledgments}
	I wish to first thank Siegfried Beckus for introducing me to this problem and suggesting the methods to tackle this question. I am also indebted to Ram Band, Siegfried Beckus and Felix Pogorezlski. Their collaboration and discussions with me on  \cite{Curr} led to this project, where many of the arguments there carried over naturally to this paper. I also wish to thank Philipp Bartmann, Daniel Lenz and  Alan Lew for discussions which contributed to this work.
	This project was supported by the Deutsche Forschungsgemeinschaft [BE 6789/1-1 to Siegfried Beckus] and ISF (grant No. 844/19).

	\newpage
	\appendix  \label{Supplement} 

	\section{Appendix}


	\subsection{Preliminaries}\label{Prelim-Sect}
	\subsubsection{Symbolic dynamical systems, dictionaries and substitutions} \label{dynam-sect}
	In this section, we  collect some basic terminology and facts regarding symbolic dynamics and substitutions, which we glossed over in the main paper. We adopt terminology from \cite[Section~II]{BBdN20} and \cite[Sections 4 and 5]{Queff10}. 
	Other recommended references on this topic are \cite[Chapter 4]{BaakeGrimm13} and \cite[Chapter 1]{Fogg}.\\
	Given an alphabet $\mathcal{A}$, the set of possible finite words over $\mathcal{A}$ is given by $\mathcal{A}^+:= \cup_{n=1}^\infty \mathcal{A}^n$. 
	We often implicitly refer to $u\in \mathcal{A}^\ell$, for some $\ell \in \mathbb{N}$, as a function $u:[0,\ell)\cap \mathbb{Z} \to \mathcal{A}$. 
	For $u_1,u_2 \in \mathcal{A}^+$, we write $u_1\prec u_2$ if $u_1$ is a subword of $u_2$.
	A \textbf{configuration} over $\mathcal{A}$, also called a two-sided infinite word, is any function $\omega: \mathbb{Z} \to \mathcal{A} $. The space of all configurations on $\mathbb{Z}$ is denoted by $\mathcal{A}^\mathbb{Z}$. We have the natural action of  $\mathbb{Z}$ on $\mathcal{A}^\mathbb{Z}$ such that, for all $m\in \mathbb{Z}$ and $\omega \in \mathcal{A}^\mathbb{Z}$, $m\cdot \omega$ is a configuration defined by
	\begin{equation}\label{Action}
		[m\cdot \omega](n):= \omega(n-m) \quad \text{for all}\quad  n\in \mathbb{Z}.
	\end{equation}
	This action is continuous with respect to the metric between configurations given by
	\begin{equation} \label{config-metric}
		d(\omega_1,\omega_2):= \inf\Big\{ \frac{1}{r+1} : \; r\geq 0 \; \text{such that} \; \omega_1 \vert_{(-r,r)}= \omega_2 \vert_{(-r,r)}    \Big \}.
	\end{equation} 
	The \textbf{orbit} of a configuration $\omega \in \mathcal{A}^\mathbb{Z}$, the closure of which is also called the hull of $\omega$, 
	is given by $Orb_\mathbb{Z}(\omega):= \{ m\cdot \omega: m\in \mathbb{Z} \}$. We call a nonempty set of configurations, $\Omega \subseteq \mathcal{A}^\mathbb{Z}$, \textbf{invariant} if $Orb_\mathbb{Z}(\omega) \subseteq \Omega$ for all $\omega \in \Omega$. An invariant nonempty subset $\Omega\subseteq \mathcal{A}^\mathbb{Z}$ is called a \textbf{subshift}, if it is also closed with respect to the metric in \eqref{config-metric}. The collections of subshifts is denoted by $\mathcal{J}$. A subshift $\Omega \in \mathcal{J}$ is called \textbf{minimal} if $\Omega =\overline{Orb_\mathbb{Z}(\omega)}$ for all $\omega\in \Omega$. A subshift $\Omega \in \mathcal{J}$ is called \textbf{periodic}, if it is minimal and finite. We note that a periodic subshift $\Omega \in \mathcal{J}$ satisfies $\Omega =Orb_\mathbb{Z}(\omega)$ for all $\omega \in \Omega$. It follows that for any periodic $\Omega \in \mathcal{J}$, there exists a $u\in \mathcal{A}^{\vert \Omega\vert}$ and $\omega' \in \Omega$, such that $\omega'\vert_{[0, \vert u\vert )}=u$ and $\omega'(n)= u(n-\lfloor \frac{n}{\vert u\vert} \rfloor)$ for all $n\in \mathbb{Z}$. We write in this case $\omega=u^\infty$. When  $\Omega \in \mathcal{J}$ is periodic, then there exists $u_\Omega\in \mathcal{A}^+$ satisfying $\Omega =Orb_\mathbb{Z}(u_\Omega^\infty)$.\\
	We now turn to discuss the notion of approximating subshifts. The space $\mathcal{J}$ is naturally endowed with the Hausdorff metric induced from $d$ given in \eqref{config-metric}. Recall that the Hausdorff metric of two nonempty sets $A,B$ is defined by
	\begin{equation} \label{Haus-metric}
		d_H(A,B):= \max \Big\{ \underset{a\in A}{\sup}\; d(a,B), \underset{b\in B}{\sup}\; d(b,A) \Big\}, \quad \text{with} \quad d(a,B):=\underset{b\in B}{\inf}d(a,b).
	\end{equation}
	The Hausdorff metric induced on $\mathcal{J}$, which we denote by $d_\mathcal{H}$, is a complete metric. We denote by $d_H$ the standard metric on subsets of $\mathbb{C}$ or $\mathbb{R}$. 
	For further details and properties of the Hausdorff distance, consult \cite{Bee93}.\\ 
	The alternative characterization of $d_\mathcal{H}$ we used in the main paper relies on the notion of dictionaries. For a configuration $\omega\in \mathcal{A}^\mathbb{Z}$ and a subshift $\Omega$, their \textbf{dictionaries} are defined accordingly as 
	\begin{equation}
		W(\omega):= \{ \omega\vert_{[n_1,n_2)}: \; n_1,n_2\in \mathbb{Z}\; n_2>n_1 \} \quad \text{and} \quad W(\Omega):=\cup_{\omega' \in \Omega} W(\omega').  
	\end{equation}
	The dictionary is comprised of the local patterns occurring in $\omega$ or $\Omega$ accordingly.
	This is also referred to as a language  in certain literature. 	 
	For every $\ell \in \mathbb{N}$, we denote $W(\Omega)_\ell:= W(\Omega)\cap \mathcal{A}^{\ell}$ and $W(\omega)_\ell:= W(\omega)\cap \mathcal{A}^\ell$. \\ 
	Subshifts and dictionaries provide equivalent definitions of essentially the same objects, with the latter often considered as simpler to handle.  
	However, to more easily obtain results from \cite{BBC19} and \cite{BecTak21}, results here are formulated using Hausdorff distance.
	These two approaches are the equivalent by either\cite[Proposition 1.3.4]{LindMarc} or \cite[Theorem 2.2]{BBdN20}. We require stricter quantitative estimates for the aforementioned correspondence, used implicitly when discussing the homeomorphism.
	To this end, we start by recalling the following Proposition, see \cite[Lemma 4.2]{BaakeGrimm13} or \cite[Proposition 5.2]{BBdN20}, 
	relating a dictionary of a configuration to a dictionary of its orbit.
	\begin{lemma}
		 \label{min-dict}
		 Let $\omega \in \mathcal{A}^\mathbb{Z}$. Then for all $m\in \mathbb{Z}$, $W(\omega)= W(m\cdot \omega)$ holds and $W(\eta)\subseteq W(\omega)$ for all $\eta \in \overline{Orb_\mathbb{Z}(\omega)}$. In particular, $W\big( \overline{Orb_\mathbb{Z}(\omega)} \big)= W(\omega)$ for all $\omega \in \mathcal{A}^\mathbb{Z}$.
	\end{lemma}
	Using Lemma \ref{min-dict}, we can conclude that 
	\begin{equation} \label{dict-ident}
		W(\Omega)= \{  \omega\vert_{[0,\ell)}:  \ell \in \mathbb{N},\; \omega \in \Omega  \} \quad \text{for any subshift} \quad \Omega \in \mathcal{J}.
	\end{equation}
	The following useful proposition follows, shown for example in \cite[Lemma 2.2]{BBC19}. 
	\begin{proposition} \label{subshift-dist-Prop}
		Let $\Omega_1,\Omega_2 \in \mathcal{J}$ and $\rho\in \mathbb{N}$. Then $d_\mathcal{H}(\Omega_1,\Omega_2) \leq\frac{1}{\rho+1}$ if and only if $W(\Omega_1)_{2\rho-1}=W(\Omega_2)_{2\rho-1}$.
	\end{proposition}
	
	\begin{remark}
		This proposition helps shed some light on the different behaviour of distance between subshifts and distance between configurations. Heuristically, the distance between two configurations corresponds to how much they agree around the origin. On the other hand, the distance between two subshifts corresponds to the lengths of words on which their dictionaries agree.
	\end{remark}
	
	We now recall relevant notions for substitutions from \cite[Section 5]{Queff10} and \cite[Section~II]{BBdN20}. A \textbf{substitution rule}, $S_0:\mathcal{A}\to \mathcal{A}^+$, is extended to maps on $\mathcal{A^+}$ and $\mathcal{A}^{\mathbb{Z}}$ by concatenation. 
	The \textbf{substitution matrix} of $S$ is the matrix $M(S)\in \mathbb{R}_{s \times s}$, such that $[M(S)]_{i,j}$ is the number of times the letter $a_i$ occurs in $S(a_j)$, where $\mathcal{A}=\{ a_1,...,a_s \}$. The matrix $M(S)$ is called the \textbf{substitution matrix} of $S$. 
	Recall that $S$ is called \textbf{primitive}, if there exist a $p\in \mathbb{N}$, such that $b\prec S^p(a)$ for all $a,b\in \mathcal{A}$. 
	The collection of legal words  is called the dictionary of $S$ and denoted $W(S)$. Equivalently,
	\begin{equation}
		 W(S)= \underset{a\in \mathcal{A}}{\bigcup}\; \underset{n\in \mathbb{N}}{\bigcup} \big\{ u\in \mathcal{A}^+:\; u\prec S^n(a) \big \}.
	\end{equation}
	For general primitive substitutions, the definition of $W(S)$ yields a unique subshift $\Omega(S)\in \mathcal{J}$, called the \textbf{substitution subshift} satisfying $W( \Omega(S) )= \Omega(S)$. 
	Using the introduced terminology and notations, the IHS approximants  are rigorously defined as $\Omega_n(\omega_0):=  \overline{ Orb_\mathbb{Z}\big( S^n(\omega_0)  \big)  }$. 
	Another classic property of $\Omega(S)$, that we rely on for the proof of Theorem \ref{1D-conv-rate}, is linear repetitivity of $\Omega(S)$. Recall that a configuration $\omega \in \mathcal{A}^\mathbb{Z}$ is called \textbf{linear repetitive} if there exists some $C\geq 1$ such that, for all $u,w\in W(S)$, if $\vert w\vert \geq C \vert u\vert $ then $u\prec w$. In this case, we say that the aforementioned $C$ is a linear repetitivity constant of $\omega$. 
	We choose some linear repetitivity constant of some $\omega \in\Omega(S)$ and denote it by $C_S\geq 1$. This $C_S$ is well defined by \cite[Lemma 2.3]{Sol98} or \cite[Proposition 6]{Dur00}. 

	\subsubsection{Directed graphs}\label{Grap-Sect}
		To clarify the terms used in Theorem \ref{1D-necess&suff}, we recall several relevant graph notions in this subsection.
		A finite set $V$, over which a graph is defined, is called a \textbf{vertex} set. A \textbf{directed graph} $G$ over $V$ is a pair $G=(V,E)$, where $E\subseteq V\times V$. The set $E$ is called the \textbf{directed edge} set. For any $v,u\in V$, we write $v\to u$, or $v\overset{G}{\to}u$, if $(v,u)\in E$. A tuple $(u_0,...,u_\ell)$ of vertices, for some $\ell \in \mathbb{N}$, is called a \textbf{directed path} if $(u_{j-1},u_j)\in E$ for all $1\leq j \leq \ell$. In that case, we say $u_0$ \textbf{leads to} $u_\ell$ and that the path \textbf{starts at} $u_0$. We say that the \textbf{length} of such a path is $\ell$. A path $(u_0,...,u_\ell)$ is said to be \textbf{closed} if $u_0=u_\ell$. \\
		Given two paths $(u_0,...,u_\ell)$ and $(v_0,...,v_{\ell'})$, we say $(v_0,...,v_{\ell'})$ is a \textbf{directed subpath} of $(u_0,...,u_\ell)$ if $\ell'\leq \ell$ and there exists some $0\leq j\leq \ell$, such that $u_{j+i}=v_i$ for all $0\leq i\leq \ell'$. We also say that $(u_0,...,u_\ell)$ contains $(v_0,...,v_{\ell'})$ as a subpath. 
		We now rigorously define $G(S)$. 
		We wish to see what are the illegal $2$-words occurring after an application of $S$ to any $u\in \mathcal{A}^2$. The following definition encodes the desired behavior.
	\begin{defn} \label{1D-ileg-graph}
		Let $S$ be a  substitution. Define a directed graph $G(S)$ with vertex set $V=\mathcal{A}^2$, such that $(u,w)$ is a directed edge in $G(S)$, denoted also by $u\to w$, if and only if $w$ is a subword of $S(u)$ and both $w,u$ are illegal words.
	\end{defn}
	We note that this definition implies that $S$ is self-correcting if and only if $G(S)$ has no closed paths.
	 \subsection{Formal statements and  proofs outline}\label{Proof-Sect}
	 \subsubsection{Convergence criteria for the scheme}
 		We outline how to show Theorem \ref{1D-necess&suff} in this subsection. 
 		As mentioned earlier,  \cite[Proposition 2.11]{BBdN20}   provides a sufficient condition for the convergence of the IHS $\Omega_n(\omega_0)$. 
 		Theorem \ref{1D-necess&suff} follows from  \cite{BBdN20} and the pigeon-hole principle. 
 		The proof of Theorem \ref{1D-necess&suff} relies on the following lemma. The following  lemma relates illegal $2$-words in $W\big( S^n(\omega_0) \big)$ to directed paths in $G(S)$.
 		\begin{lemma} \label{path-corr}
 			Let $S$ be a primitive aperiodic one-dimensional substitution over an alphabet $\mathcal{A}$ and let $\omega_0 \in \mathcal{A}^\mathbb{Z}$.
 			\begin{enumerate}
 				\item If $u\in W\big( S(\omega_0) \big)_2\setminus W(S)$, then there exists some $w\in W(\omega_0)_2\setminus W(S)$ such that $u\preceq S(w)$.
 				\item If $u\in W\big( S^n(\omega_0) \big)_2\setminus W(S)$, then there exist $\{ u_j \}_{j=0}^{n-1}\subseteq \mathcal{A}^2$ such that
 				\[ u_j \in W\big( S^j(\omega_0) \big)_2\setminus W(S) \quad \text{and} \quad u_{j+1} \prec S(u_j)  \quad \text{for all} \quad 0\leq j \leq n-1,   \]
 				where $u_n:=u$.
 				In particular, there exists some $u_0 \in W(\omega_0)_2\setminus W(S)$ and a path of length $n$ in $G(S)$, $(u_0,...,u_{n-1},u)$,  such that $u_j\to u_{j+1}$ for all $0\leq j\leq n-1$. 
 			\end{enumerate} 
 		\end{lemma}
 		The proofs of the Lemma \ref{path-corr} and Theorem \ref{1D-necess&suff} are similar to the proofs of analogous claims in \cite{Curr}, and so we omit them. 
 		\subsubsection{Convergence rate} \label{Rate-Sect}
 		We outline how to prove the estimates in \eqref{Conv-rate} in this subsection. We start by restating the estimates in \eqref{Conv-rate} more precisely. 
 		\begin{Theorem}\label{1D-conv-rate}
 			Let $S$ be a primitive one-dimensional substitution over an alphabet $\mathcal{A}$, with PF-eigenvalue $\theta_S$.  Given an initial configuration $\omega_0\in \mathcal{A}^\mathbb{Z}$,  consider the generated IHS $\Omega_n(\omega_0)$. 
 			Then there exists a constant $C_1=C_1(S)>0$ such that if $\Omega_n(\omega_0)$ are good approximating IHS, then
 			\[  d_\mathcal{H}\big( \Omega_n(\omega_0), \Omega(S) \big) \leq \frac{C_1}{\theta_S^n} \quad \text{for all} \quad n\in \mathbb{N}. \]
 			If furthermore, $\omega_0$ is a periodic configuration, then there exists a constant $C_2=C_2(\omega_0;S)>0$ such that
 			\[ \frac{C_2}{\theta_S^n}   \leq  d_\mathcal{H}\big( \Omega_n(\omega_0), \Omega(S) \big) . \]
 		\end{Theorem}  
 		We once again only outline the proof for Theorem \ref{1D-conv-rate}, due to similarities with arguments in \cite{Curr}. 
 		Theorem \ref{1D-conv-rate} is proved using the following lemmas. 		
 		\begin{lemma} \label{upp-bound}
 			Let $S$ be a primitive one-dimensional substitution over an alphabet $\mathcal{A}$, with linear repetitive constant $C_S$ and PF-eigenvalue $\theta_S$. Let $r\in \mathbb{N}$, $\omega_0 \in \mathcal{A}^\mathbb{Z}$ and consider the generated IHS $\Omega_n(\omega_0)$.
 			\begin{enumerate}
 				\item If $n\geq N_1(r)$, then $W\big( \Omega_n(\omega_0) \big)_r\supseteq W(S)_r$, where $N_1(r):=\frac{\log r}{\log \theta_S}+ \frac{\log C_S -\log \check{C}(S)}{\log \theta_S}$.
 				\item If the IHS $\Omega_n(\omega_0)$ are good approximants and $n>N_2(r)$, then $W\big( \Omega_n(\omega_0) \big)_r\subseteq W(S)_r$, where $N_2(r):=\frac{\log r}{\log \theta_S}- \frac{\log (2\check{C}(S) )}{\log \theta_S}+\vert \mathcal{A}\vert^2$.
 			\end{enumerate}
 		\end{lemma}
 		This lemma provides us with finer quantitative bounds on $d_{\mathcal{H}}\big( \Omega_n(\omega_0), \Omega(S) \big)$.
 		By Proposition \ref{subshift-dist-Prop}, we can conclude that $c(\Omega_1,2r-1)=c(\Omega_2,2r-1)$, whenever $\Omega_1,\Omega_2\in \mathcal{J}$ satisfy $d_\mathcal{H}(\Omega_1,\Omega_2) \leq\frac{1}{r+1}$.
 		The lower bounds on $d_{\mathcal{H}}\big( \Omega_n(\omega_0), \Omega(S) \big)$ for periodic $\omega_0$ comes from the complexity function. 
 		For a subshift $\Omega\in \mathcal{J}$, its \textbf{complexity function} is defined to be
 		\begin{equation} \label{complex-func}
 			c(\Omega,\cdot):\mathbb{N} \to \mathbb{N}, \quad c(\Omega, r):= \big\vert W(\Omega)_r \big\vert.
 		\end{equation}
 		
 		\begin{lemma}\label{complex-bound}
 			Let $S$ be a primitive substitution over an alphabet $\mathcal{A}$, with PF-eigenvalue $\theta_S$. Let $\omega_0\in \mathcal{A}^\mathbb{Z}$ be a periodic initial configuration, with $\omega_0=u^\infty$ for $u\in \mathcal{A}^+$. Consider the generated IHS $\Omega_n(\omega_0)$. Then $c\big( \Omega_n(\omega_0),r \big)\leq \hat{C}(S)\cdot \theta_S^n \cdot \vert u\vert$ holds for all $r\in \mathbb{N}$.
 		\end{lemma}
 		Combining Lemma \ref{complex-bound} with the Coven--Hedlund lower bound \cite[Corollary 2.12]{Coven73} for complexities of aperiodic systems, one can derive the lower bound for $d_{\mathcal{H}}\big( \Omega_n(\omega_0), \Omega(S) \big)$.
 		\subsubsection{Specrtal convergence rates}
 		For a  Schrödinger operator $H_\omega$, as in \eqref{Schro-defn}, we denote the arising spectrum by $\sigma(H_\omega)$.
 		A precise version of the results in \cite{BBC19} and \cite{BecTak21}, is given in the following proposition.
 		\begin{proposition}\cite[Proposition 1.1]{BBC19} \label{Lipschitz-Bound}
 			Let $H_\omega$ be a Schrödinger operator as in \eqref{Schro-defn}. Then there exists a constant $C>0$ such that
 			\begin{equation}
 				d_H \big( \sigma(H_{\omega_1}), \sigma(H_{\omega_2})  \big) \leq C d_\mathcal{H} \Big( \overline{Orb_\mathbb{Z}(\omega_1)}, \overline{Orb_\mathbb{Z}(\omega_2)}  \Big) \quad \text{for any} \quad \omega_1,\omega_2\in \mathcal{A}^\mathbb{Z}.
 			\end{equation}
 		\end{proposition} 
 		Finally, we precisely state the exponentially fast approximation result for a self-correcting substitution.
 		\begin{proposition} \label{Approx-Method}
 			Let $S$ be a self-correcting primitive substitution over an alphabet $\mathcal{A}$, with PF-eigenvalue $\theta_S$. Then there exists a constant $\tilde{C}>0$ such that
 			\[ d_H \Big( \sigma( H_{S^n(\omega_0)} ), \sigma ( H_{\omega_S} )  \Big) \leq \frac{\tilde{C}}{ \theta_S^n  } \quad \text{for any} \quad \omega_0\in \mathcal{A}^\mathbb{Z}, \]
 			where $\omega_S\in \mathcal{A}^{\mathbb{Z}}$ is a configuration with substitution symmetry. 
 		\end{proposition}
 		This result follows directly by combining \cite[Proposition 1.1]{BBC19} and Theorem \ref{1D-conv-rate}, so we leave its proof to the reader.
 		
 	\subsection{Further examples of substitutions}
\begin{figure}[H]
 			\vspace*{2cm}
 			\begin{center}
 				\includegraphics[width=0.5\linewidth]{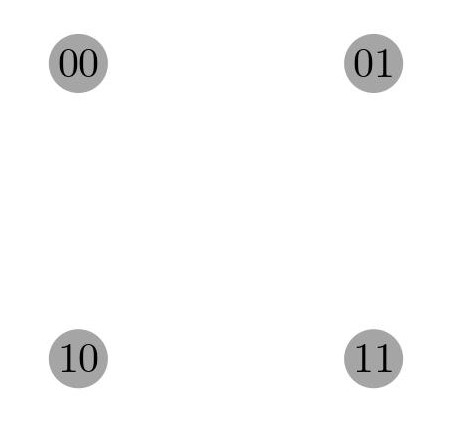}
 				\caption{$G(S_{TM})$ for the Thue--Morse substitution.}
 				\label{Thue-Morse-graph} 
 			\end{center}
 			
 		\end{figure}

\begin{figure}[H]
 			\begin{center}
 				\includegraphics[width=1\linewidth]{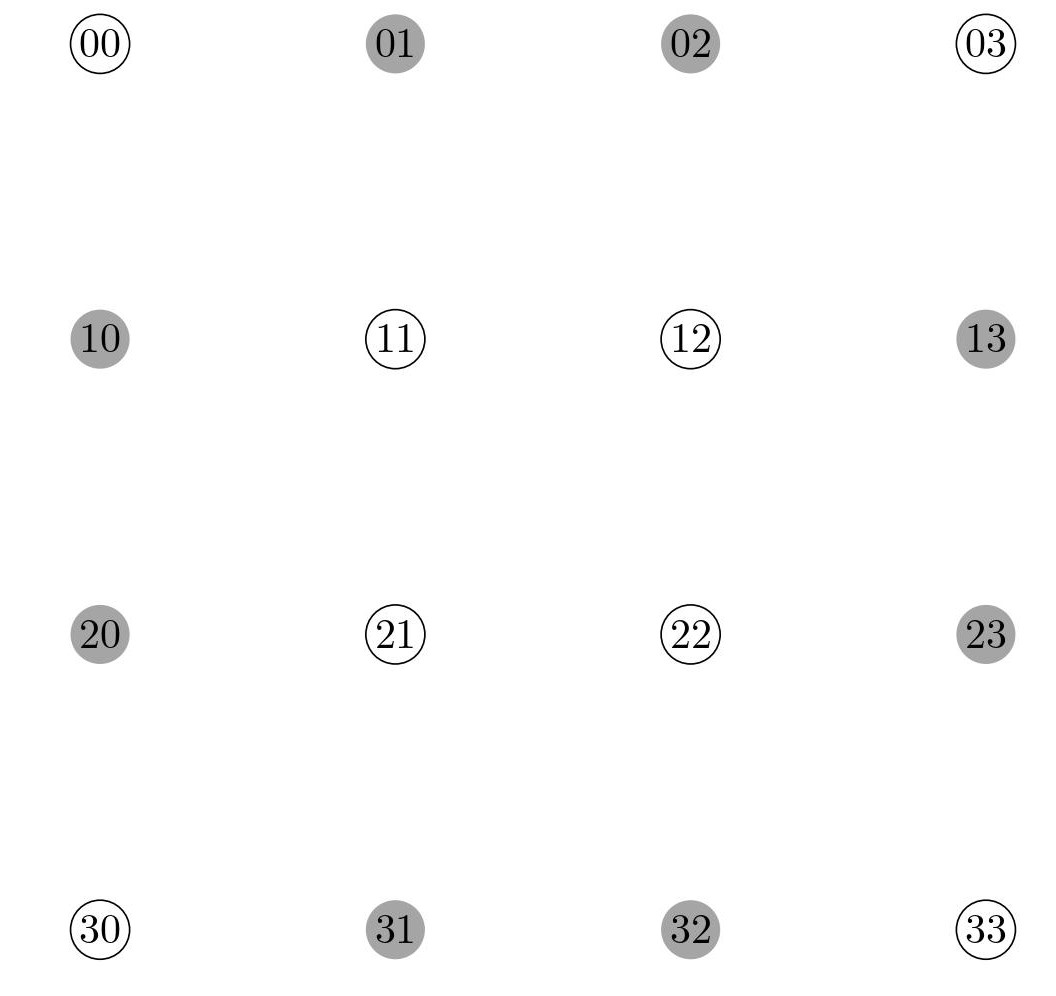}
 				\caption{\centering $G(S_{GRS})$ for the Golay--Rudin--Shapiro substitution.} 
 				\label{Gol-Rud-graph}
 			\end{center}
 		\end{figure}

 	In this subsection we consider some more examples of substitutions and review their resulting $G(S)$ graphs. One can see from the lack edges in the following figures, that all these substitutions are indeed self-correcting.
 	\begin{example}\label{Period doubling}
 		Let $\mathcal{A}=\{0,1\}$. The period doubling substitution rule is given by
 		\[ 0 \mapsto 01\quad \text{and} \quad 1 \mapsto 00. \]
 		Its $G(S_{PD})$ graph is the same as that for the Fibonacci substitution and is given in Figure \ref{Fibo-graph}.
 	\end{example}
 	\begin{example}\label{Thue-Morse}
 		Let $\mathcal{A}=\{0,1\}$. The Thue--Morse substitution rule is given by
 		\[ 0 \mapsto 01\quad \text{and} \quad 1\mapsto 10. \]
 		Its legal $2$-words are all $\mathcal{A}^2$ and its $G(S_{TM})$ graph is given in Figure \ref{Thue-Morse-graph}.
 	\end{example}
 	\begin{example} \label{Gol-Rud-Shap}
 		Let $\mathcal{A}=\{0,1,2,3\}$. The Golay--Rudin--Shapiro substitution rule is given by
 		\[ 0\mapsto 01, \quad 1 \mapsto 02,\quad 2 \mapsto 31 \quad  \text{and} \quad  3 \mapsto  32. \]
 		Its $G(S_{GRS})$ graph is given in Figure \ref{Gol-Rud-graph} with the legal $2$-words colored again in gray. 
 	\end{example}
 	\begin{multicols}{2}

 		\columnbreak

 	\end{multicols}

	\bibliography{1DProceedRev}
	\bibliographystyle{alpha}
	
\end{document}